\documentclass[12pt]{article}
\usepackage{amssymb, amsmath, amsfonts,dsfont, mathrsfs,euscript,eufrak,colonequals}
\usepackage[T1]{fontenc}
\usepackage{cmap}
\usepackage{textcomp}
\usepackage{latexsym}
\usepackage{indentfirst}

\newtheorem{Theorem}{Theorem}

\newcommand{\N}{\mathbb N}
\newcommand{\R}{\mathbb R}
\newcommand{\Z}{\mathbb Z}

\newcommand{\defeq}{\colonequals}
\newcommand{\dd}{d}

\newcommand{\la}{\langle}
\newcommand{\ra}{\rangle}

\begin{document}
	\author{A. A. Khartov\footnote{Smolensk State University, 4 Przhevalsky st., 214000 Smolensk, Russia, e-mail: \texttt{alexeykhartov@gmail.com} }}

	\title{A criterion of quasi-infinite divisibility for discrete laws}
	
\maketitle

\begin{abstract}
	We consider arbitrary discrete probability laws on the real line. We obtain a criterion of their belonging to a new class of quasi-infinitely divisible laws, which is a wide natural extension of the class of well known infinitely divisible laws through the L\'evy type representations.
\end{abstract}

\textit{Keywords and phrases}: discrete probability laws, characteristic functions, spectral L\'evy type representations, quasi-infinitely divisible laws.

\section{Introduction and problem setting}

This note is devoted to the question posed in \cite{AlexeevKhartov} about the description of quasi-infinitely divisibile laws within the class of arbitrary univariate  discrete probability laws. 

The class of quasi-infinitely divisible laws is a rather wide extention of the class of well known infinitely divisible laws.  Let $F$ be a distribution function on the real line with the characteristic function $f$. Following \cite{LindPanSato}, $F$ and $f$ (and the corresponding probability law)  are called \textit{quasi-infinitely divisible} if $f$ admits the following L\'evy type representation:
\begin{eqnarray}\label{repr_f}
	f(t) = \exp\Biggr\{it\gamma - \frac{\sigma^2 t^2}{2} + \int\limits_{\R\setminus\{0\}} \bigl(e^{itu}-1- it\sin  u\bigr)\dd \Lambda(u)\Biggr\}, \quad t\in\R,
\end{eqnarray} 
for some $\gamma\in\R$ (the set of real numbers), $\sigma^2 \geqslant 0$, and function $\Lambda$ (the \textit{L\'evy spectral function}), which has a finite total variation on every interval $(-\infty,-r]$ and $[r,\infty)$, $r>0$, and it satisfies
\begin{eqnarray*}
	\lim_{u \to -\infty}\Lambda(u) = \lim_{u\to+\infty}\Lambda(u)  = 0,\quad\text{and}\quad \int\limits_{0<|u|<\delta} x^2 \dd|\Lambda|(u)<+\infty,\quad\delta>0.
\end{eqnarray*}
Here $f$ uniquely determines the triplet $(\gamma, \sigma^2, \Lambda)$. Quasi-infinitely divisible  $F$ and $f$ are \textit{rationally infinitely divisible}, i.e. there exist infinitely divisible distribution functions $F_1$ and $F_2$ with characteristic functions $f_1$ and $f_2$, respectively, such that $F_1=F*F_2$, where ``$*$'' is the convolution, or in the equivalent form: $f(t)=f_1(t)/f_2(t)$, $t\in\R$. It is clear that every infinitely divisible distribution function is quasi-infinitely divisible with non-decreasing $\Lambda$ on every interval $(-\infty, 0)$ and $(0,\infty)$. The opposite is not true, the examples can be found in the classical monographs \cite{GnedKolm}, \cite{LinnikOstr}  and \cite{Lukacs}.

The first detailed analysis of quasi-infinitely divisible laws based on their L\'evy type representations was performed in \cite{LindPanSato}, where  the authors studied  questions concerning  supports, moments, continuity, and weak convergence for these laws. These results were generalized and complemented in the papers  \cite{AlexeevKhartov}, \cite{AlexeevKhartov2}, \cite{Berger}, \cite{BergLind2}, \cite{BergLind},  \cite{Khartov}, \cite{KhartAlexeev} and \cite{Kutlu}. The quasi-infinitely divisible laws have already found a lot of interesting applications  (see \cite{BergLind}, \cite{ChhDemniMou}, \cite{LindSato}, \cite{Nakamura}, \cite{ZhangLiuLi} and the references given there).

It is rather natural and interesting to investigate criteria of quasi-infinite divisibility for probability laws. The most deep results here  were obtained  for univariate discrete  laws  in \cite{AlexeevKhartov}, \cite{KhartAlexeev}, and \cite{LindPanSato}.  In particular, following \cite{LindPanSato}, a discrete lattice probability law  is quasi-infinitely divisible if its characteristic function $f$ has no zeroes on the real line, i.e. $f(t)\ne 0$, $t\in\R$. For arbitrary discrete laws only the sufficient condition is known due to  \cite{AlexeevKhartov}. Namely, if $f$ is separated from zero, i.e. $|f(t)|\geqslant\mu>0$ for some $\mu$ and all $t\in\R$, then the law is quasi-infinitely divisible (see the next section for more details). The latter result generalizes the previous, because in the lattice case the function $t\mapsto |f(t)|$, $t\in\R$, is continuous and periodic, and the absence of zeroes is equivelent to the separatness from zero over the period segment (see \cite{AlexeevKhartov} for more details). Anyway, however, the following interesting question remains here (it was posed in \cite{AlexeevKhartov}). Are all characteristic functions of quasi-infinitely divisible laws separated from zero? In this note we will show that it is true. So we obtain a criterion of quasi-infinite divisibility for arbitrary discrete laws.

We will use the following notation. We denote by $\R$, $\Z$ and $\N$ the sets of real numbers, integers, and positive integers, respectively. We write $a_n\sim b_n$ for number sequences $a_n$ and $b_n$, $n\in\N$,  if $a_n/b_n\to 1$, $n\to\infty$.

\section{The result}
Let us consider an arbitrary discrete  distribution function
\begin{eqnarray*}
	F(x) \defeq \sum_{k\in\N:\, x_k\leqslant x}p_{x_k}, \quad x \in \R,
\end{eqnarray*} 
where  $x_k$, $k\in\N$, are distinct real numbers, $p_{x_k}\geqslant 0$, $k\in \N$, and $\sum_{k=1}^\infty p_{x_k}=1$. Let $f$ be the characteristic function of $F$:
\begin{eqnarray}\label{def_f}
	f(t)=\sum_{k\in\N} p_{x_k} e^{it x_k},\quad t\in\R.
\end{eqnarray}
Let $X$ be the set of all points of growth of $F$, i.e. $X\defeq \{x_k:  p_{x_k}>0,\, k\in\N \}\ne \varnothing$. Let us introduce the set of all finite $\Z$--linear combinations of elements from the set $X$:
\begin{eqnarray*}
	\la X\ra\defeq \biggl\{\sum_{k=1}^m c_k x_k:\,   c_k \in \Z,\, x_k\in X,\, m\in\N\biggr\}.
\end{eqnarray*}
In other words, $\la X\ra$ is a module over the ring $\Z$ with the generating set $X$. Since $X\ne \varnothing$,  $\la X\ra$ is an infinite countable set. The following result was obtained in \cite{AlexeevKhartov}.

\begin{Theorem}\label{th_Repr}
	Suppose that  $\inf_{t\in\R} |f(t)|>0$, i.e. there exists $\mu>0$ such that $|f(t)|\geqslant \mu>0$ for all $t\in\R$. Then $f$ admits the following representation
	\begin{eqnarray}\label{th_Repr_eq}
		f(t) = \exp\biggr\{ it\gamma_0 + \sum_{u \in \la X\ra\setminus \{0\}}\lambda_{u}(e^{itu}-1)\biggr \}, \quad t\in\R,
	\end{eqnarray}
	where $\gamma_0 \in \la X\ra$, $\lambda_u \in \R$ for all $u\in\la X\ra\setminus \{0\}$, and $\sum_{u \in \la X\ra\setminus \{0\}}|\lambda_u| < \infty$. Here $F$ is quasi-infinitely divisible and $f$ has representation \eqref{repr_f} with	
	\begin{eqnarray*}
		\gamma= \gamma_0+\!\!\!\sum_{u \in \la X\ra\setminus \{0\}}\lambda_{u}\sin(u),\qquad \sigma^2=0,\qquad \Lambda(x)=
		\begin{cases}
			\sum\limits_{u\in\la X\ra:\, u\leqslant x} \lambda_{u},&\quad x<0,\\
			-\!\!\!\!\!\!\sum\limits_{u\in\la X\ra:\, u > x} \lambda_{u},&\quad x>0.
		\end{cases}
	\end{eqnarray*}	
\end{Theorem}

Thus the condition $\inf_{t\in\R} |f(t)|>0$ implies the  quasi-infinite  divisibility  of $F$.  We complement this fact by the following proposition.

\begin{Theorem}\label{th_inff0}
	Suppose that $\inf_{t\in\R} |f(t)|=0$. Then $F$ is not quasi-infinitely divisible.
\end{Theorem}
\textbf{Proof.}  It is known that characteristic functions of quasi-infinitely divisible laws have no zeroes on  the real line (see \cite{LindPanSato}). Therefore if there exists $t_0\in\R$ such that $f(t_0)=0$, then $F$ is not quasi-infinitely divisible. Hence we focus only on the case $f(t)\ne 0$ for all $t\in\R$ with $\inf_{t\in\R} |f(t)|=0$.
 
Suppose, contrary to our claim, that $F$ is a quasi-infinitely divisible distribution function. Then $f$ admits  representation \eqref{repr_f} with some triplet $(\gamma,\sigma^2, \Lambda)$. Let us consider for fixed $\tau\in\R$ the following function
\begin{eqnarray*}
	\psi_\tau (t)\defeq \dfrac{f(t+\tau)f(t-\tau)}{f(t)^2}=\exp\biggl\{  -\sigma^2\tau^2 +  2\int\limits_{\R\setminus\{0\}} e^{itx} \bigr(\cos(\tau x)-1 \bigr) \dd \Lambda(x) \biggr\},\quad t\in\R.
\end{eqnarray*}
Due to the conditions for $\Lambda$, for any fixed $\tau\in\R$ we have 
\begin{eqnarray}\label{conc_suppsi}
	\sup_{t\in\R}|\psi_\tau (t)|\leqslant \exp\biggl\{ -\sigma^2\tau^2 +  2\int\limits_{\R\setminus\{0\}} \bigr(1-\cos(\tau x)\bigr) \dd |\Lambda|(x) \biggr\} <\infty.
\end{eqnarray}
 
Let $(t_m)_{m\in\N}$ be an increasing sequence such that $t_m\to+\infty$ and $f(t_m)\to 0$, $m\to\infty$. Let us define $\varphi_m(\tau)\defeq f(t_m+\tau)$, $\tau\in\R$, $m\in\N$. Since $f$ is an almost periodic function by \eqref{def_f}, the sequence $(\varphi_m)_{m\in\N}$ is relatively compact in the topology of uniform convergence on $\R$ (see \cite{Levitan} pp. 23--24). So there exists a subsequence $(\varphi_{m_l})_{l\in\N}$ that  uniformly converges to an almost periodic function $\varphi$, i.e.
\begin{eqnarray}\label{conc_conv}
	\sup_{\tau\in\R}|f(t_{m_l}+\tau)-\varphi(\tau)|\to 0,\quad l\to\infty.
\end{eqnarray}
From this we have
\begin{eqnarray}\label{conc_conv_pmphi}
	\sup_{\tau\in\R}|f(t_{m_l}+\tau)f(t_{m_l}-\tau)-\varphi(\tau)\varphi(-\tau)|\to 0,\quad l\to\infty.
\end{eqnarray}
Indeed, for any $\tau\in\R$ and $l\in\N$ it holds that
\begin{eqnarray*}
	f(t_{m_l}+\tau)f(t_{m_l}-\tau)-\varphi(\tau)\varphi(-\tau)= \bigl(f(t_{m_l}+\tau)-\varphi(\tau)\bigr)f(t_{m_l}-\tau)+\bigl(f(t_{m_l}-\tau)-\varphi(-\tau)\bigr)\varphi(\tau).
\end{eqnarray*}
Due to $|f(t_{m_l}-\tau)|\leqslant 1$ and $|\varphi(\tau)|\leqslant 1$ for any $\tau\in\R$ (see \eqref{def_f} and \eqref{conc_conv}), we get the inequality
\begin{eqnarray*}
	\bigl|f(t_{m_l}+\tau)f(t_{m_l}-\tau)-\varphi(\tau)\varphi(-\tau)\bigr|&\leqslant&  \bigl|f(t_{m_l}+\tau)-\varphi(\tau)\bigr|+\bigl|f(t_{m_l}-\tau)-\varphi(-\tau)\bigr|\\
	&\leqslant& 2\sup_{\tau\in\R}|f(t_{m_l}+\tau)-\varphi(\tau)|,
\end{eqnarray*}
which (together with \eqref{conc_conv}) implies \eqref{conc_conv_pmphi}.

Suppose that  $\varphi(\tau)\varphi(-\tau)= 0$ for any $\tau\in\R$, i.e.
\begin{eqnarray*}
	\sup_{\tau\in\R}|f(t_{m_l}+\tau)f(t_{m_l}-\tau)|\to 0,\quad l\to\infty.
\end{eqnarray*}
In particular, for any fixed $s\in\R$  we have
\begin{eqnarray*}
	\bigl|f(t_{m_l}+\tau)f(t_{m_l}-\tau)\bigr|\biggl|_{\tau=t_{m_l}+s}= |f(2t_{m_l}+s)f(-s)|\to 0,\quad l\to\infty.
\end{eqnarray*}
Since $f(t)\ne 0$, $t\in\R$, we conclude that $|f(2t_{m_l}+s)|\to 0$, $l\to\infty$ for any fixed $s\in\R$.  Next, from the sequence $(m_l)_{l\in\N}$ we can choose a  subsequence $(m'_l)_{l\in\N}$ such that the functions $s\mapsto f(2t_{m'_l}+s)$ uniformly converge on $\R$ as $l\to\infty$ (see \cite{Levitan} pp. 23--24). It is clear that their uniform limit is zero, i.e. we have
\begin{eqnarray*}
	\sup_{s\in\R}|f(2t_{m'_l}+s)|\to 0,\quad l\to\infty.
\end{eqnarray*}
So, in particular,  $|f(0)|=|f(2t_{m'_l}-2t_{m'_l})|\to 0$, $l\to\infty$, but $f(0)=1$, a contradiction.  

Thus there exists $\tau\in\R$ such that $\varphi(\tau)\varphi(-\tau)\ne 0$. Then, due to \eqref{conc_conv_pmphi} and the convergence $f(t_{m_l})\to 0$, $l\to\infty$, for this $\tau$ we have
\begin{eqnarray*}
	|\psi_{\tau} (t_{m_l})|= \dfrac{|f(t_{m_l}+\tau)f(t_{m_l}-\tau)|}{|f(t_{m_l})|^2}  \sim \dfrac{|\varphi(\tau)\varphi(-\tau)|}{|f(t_{m_l})|^2}\to\infty,\quad l\to\infty.
\end{eqnarray*}
Thus we found $\tau$ such that $\sup_{t\in\R} |\psi_{\tau}(t)|=\infty$. This contradicts \eqref{conc_suppsi}.  Therefore $f$ can not be the characteristic function of  a quasi-infinitely divisible law.\quad $\Box$\\

The examples of $f$ satisfying $\inf_{t\in\R}|f(t)|=0$ were considered in \cite{AlexeevKhartov} and \cite{KhartAlexeev}.

Thus Theorem \ref{th_Repr} and Theorem \ref{th_inff0} yield the following criterion: \textit{discrete distribution function $F$ is quasi-infinitely divisible if and only if $\inf_{t\in\R} |f(t)|>0$, i.e. there exists $\mu>0$ such that $|f(t)|\geqslant \mu>0$ for all $t\in\R$.} Moreover, now we know that the characteristic functions of all discrete quasi-infinitely divisible laws admit representation \eqref{th_Repr_eq} with some $\gamma_0$ and $\lambda_u$.

\section{Acknowledgments}
The author is supported by the Ministry of Science and Higher Education of the Russian Federation, agreement 075-15-2019-1620 date 08/11/2019.

\end{document}